\documentclass[11pt,a4paper,reqno]{amsart} 

\usepackage{amssymb,graphicx}

\usepackage{amssymb}

\newcommand{\koniec}{\begin{flushright}  $\Box $ \end{flushright}}
\newtheorem{theo}{Theorem}[section] 
\newtheorem{prop}[theo]{Proposition}  
\newtheorem{lemma}[theo]{Lemma}
\newtheorem{defi}[theo]{Definition}

\newcounter{mnotecount}[section]

\renewcommand{\themnotecount}{\thesection.\arabic{mnotecount}}

\newcommand{\mnote}[1]
{\protect{\stepcounter{mnotecount}}$^{\mbox{\footnotesize
$
\bullet$\themnotecount}}$ \marginpar{
\raggedright\tiny\em
$\!\!\!\!\!\!\,\bullet$\themnotecount: #1} }

\newcommand{\C}{\mathbb{C}}
\newcommand{\PP}{\mathbb{P}}

\def\p{\partial}
\def\be{\begin{equation}}
\def\ee{\end{equation}}

\def\bea{\begin{eqnarray}}
\def\eea{\end{eqnarray}}

\begin{document}
\date{January 24th, 2018}
\title{Metrisability of Painlev\'e equations}
\author{}
\author{Felipe Contatto }
\author{Maciej Dunajski}
\address{Department of Applied Mathematics and Theoretical Physics\\ 
University of Cambridge\\ Wilberforce Road, Cambridge CB3 0WA\\ UK.}
\email{felipe.contatto@damtp.cam.ac.uk, m.dunajski@damtp.cam.ac.uk}
\begin{abstract} 
We solve the metrisability problem for the six  Painlev\'e equations, 
and more generally for all 2nd order ODEs with Painlev\'e property,
and determine for which of these equations their integral curves are geodesics
of a (pseudo) Riemannian metric on a surface. 
\end{abstract}   
\maketitle
\section{Introduction}
A geometric approach to nonlinear 2nd order ODEs was 
initiated in the works of Liouville \cite{Liouville1886} and developed  by 
Cartan.
A general 2nd order ODE defines a path geometry on a surface $U$ coordinatised by the dependent and independent variables: there is a unique integral curve
through each point of $U$ in each direction.
The paths are unparametrised geodesics of a torsion--free connection $\nabla$ on $TU$ with Christoffel symbols $\Gamma_{ab}^c$ if and only if the ODE is of the form
\begin{equation}
\label{ODE2}
\frac{d^2 y}{d x^2}=A_3(x, y)
\Big(\frac{d y}{d x}\Big)^3
+A_2(x, y)\Big(\frac{d y}{d x}\Big)^2
+A_1(x, y)\Big(\frac{d y}{d x}\Big)+A_0(x, y)
\end{equation}
where
\be
\label{expressions_forA}
A_0=-\Gamma^2_{11},\quad A_1=\Gamma^1_{11}-2\Gamma^2_{12}, \quad 
A_2=2\Gamma^1_{12}-\Gamma^2_{22},\quad A_3=\Gamma^1_{22}.
\ee
Conversely, with any ODE of the form 
(\ref{ODE2}) one can associate 
a {\em projective structure} \cite{C22, Thomas1925}, that is an equivalence class of
torsion--free connections which share the same unparametrised geodesics. Two connections
$\nabla$ and $\hat{\nabla}$
belong to the same projective equivalence class if 
their geodesic flows on $TU$ project to the same foliation of $\PP(TU)$. Equivalently, there exists a 
one-form $\Upsilon$ on $U$ such that 
\be
\label{projequiv}
\hat\Gamma^a_{bc}=\Gamma^a_{bc}+\Upsilon_b\delta^a_c+\Upsilon_c\delta^a_b.
\ee
\begin{defi}
A second order ODE is called metrisable if its integral curves are unparametrised 
geodesics of a Levi--Civita connection of some (pseudo) Riemannian metric.
\end{defi}
A problem of characterising metrisable ODEs by differential invariants was posed
by Roger Liouville \cite{Liouville1886}, 
who has reduced it to an overdetermined system of linear PDEs (see Theorem
\ref{theo_liouville} in the next section). The complete solution was provided
relatively recently \cite{BDE}, where it was
shown that an ODE is metrisable if and only
if three point invariants of differential 
orders five and six vanish, and certain genericity assumptions hold.

\vskip5pt
A different approach was
developed by Painlev\'e, Kowalevskaya and Gambier
who studied 2nd order ODEs in the complex domain \cite{painleve, Ince}.
\begin{defi}
The ODE $y''=R(x, y, y')$, where $R$ is a rational function of $y$ and $y'$ has the Painlev\'e property (PP) if its movable
singularities (i.e. singularities whose locations depend on the initial conditions) are  poles.
\end{defi}
The solutions of equations with Painlev\'e property
are single--valued thus giving rise to proper
functions on $\C$. There exists fifty canonical types of second
order ODEs with PP up to the change of variables
\be
\label{homot}
y\rightarrow Y(x, y)=\frac{a(x)y+b(x)}{c(x)y+d(x)}, \quad
x\rightarrow X(x)=\phi(x),
\ee
where functions $(a, b, c, d, \phi)$ are analytic in $x$. Forty--four of these are solvable in terms of `known'
functions (sine, cosine, elliptic functions
or in general solutions to linear ODEs). 
The remaining six  types define new transcendental functions, and 
are given by the Painlev\'e equations  
\begin{align*}
y''&=6y^2+x &\text{PI}&\\
y''&=2y^3+xy+\alpha &\text{PII}&\\
y''&=\frac{1}{y}y'^2-\frac{1}{x}y'+\alpha \frac{y^2}{x}+\frac{\beta}{x}+\gamma y^3+\frac{\delta}{y} &\text{PIII}&\\
y''&=\frac{1}{2y}y'^2+\frac{3}{2}y^3+4xy^2+2(x^2-\alpha)y+\frac{\beta}{y} &\text{PIV}&\\
y''&=\left(\frac{1}{2y}+\frac{1}{y-1}\right)y'^2-\frac{1}{x}y'+\frac{(y-1)^2}{x^2}\left(\alpha y+\frac{\beta}{y}\right)+\gamma\frac{y}{x}+\delta \frac{y(y+1)}{y-1} &\text{PV}& \\
y''&=\frac{1}{2}\left(\frac{1}{y}+\frac{1}{y-1}+\frac{1}{y-x} \right)y'^2-\left(\frac{1}{x}+\frac{1}{x-1}+\frac{1}{y-x}\right)y'+& & \\
&+\frac{y(y-1)(y-x)}{x^2(x-1)^2}\left[\alpha+\beta\frac{x}{y^2}+\gamma\frac{x-1}{(y-1)^2}+\delta\frac{ x(x-1)}{(y-x)^2} \right] &\text{PVI}.&
\end{align*}
Here $\alpha, \beta, \gamma, \delta$ are constants. Thus PVI belongs to a four--parameter family of ODEs, etc. Some work towards characterising the 
Painlev\'e equations by point invariants of (\ref{ODE2}) has been done in
\cite{kamran, HietDry2002,kartak}.
\vskip5pt 
The aim of this paper is to determine which of the Painlev\'e equations are metrisable. 
In the next Section we shall prove the following 
\begin{theo}
\label{main_theo}
The only metrisable Painlev\'e equations are
\begin{enumerate}
\item Painlev\'e $III$, where $\alpha=\gamma=0$ or
$\beta=\delta=0$.
\item Painlev\'e $V$,  where
$\gamma=\delta=0$.
\item Painlev\'e $VI$,  where
$\alpha=\beta=\gamma=0$ and $\delta=1/2$.
\end{enumerate} 
If $\alpha=\beta=\gamma=\delta=0$ then the projective
structures defined by PIII and PV are flat. The
metrisable PVI projective structure is also flat.
\end{theo}
The flatness of a projective structure is equivalent
to the existence of a point transformation
$(x, y)\rightarrow (X(x, y), Y(x, y))$ such that the corresponding ODE (\ref{ODE2}) 
becomes\footnote{A second order ODE $y''=R(x,y,y')$ is equivalent to $(\ref{flat_ode})$ under a point transformation, if and only if it is of the form (\ref{ODE2}) and the following quantities, called \textit{Liouville invariants} 
vanish 
$$
L_1=\frac{2}{3}\dfrac{\partial^2A_1}{\partial x\partial y}-\frac{1}{3}\dfrac{\partial^2A_2}{\partial x^2}-\dfrac{\partial^2A_0}{\partial y^2}+A_0\dfrac{\partial A_2}{\partial y}+A_2\dfrac{\partial A_0}{\partial y}-A_3\dfrac{\partial A_0}{\partial x}-2 A_0 \dfrac{\partial A_3}{\partial x}-\frac{2}{3}A_1\dfrac{\partial A_1}{\partial y}+\frac{1}{3}A_1\dfrac{\partial A_2}{\partial x},
$$
$$
L_2=\frac{2}{3}\dfrac{\partial^2A_2}{\partial x\partial y}-\frac{1}{3}\dfrac{\partial^2A_1}{\partial x^2}-\dfrac{\partial^2A_3}{\partial x^2}-A_3\dfrac{\partial A_1}{\partial x}-A_1\dfrac{\partial A_3}{\partial x}+A_0\dfrac{\partial A_3}{\partial y}+2 A_3 \dfrac{\partial A_0}{\partial y}+\frac{2}{3}A_2\dfrac{\partial A_2}{\partial x}-\frac{1}{3}A_2\dfrac{\partial A_1}{\partial y}.
$$}
\be
\label{flat_ode}
\frac{d^2 Y}{d X^2}=0.
\ee 
In Section \ref{reducibility} we shall clarify a 
connection between the metrisability
of Painlev\'e equations and the existence of first
integrals: all metrisable cases are reducible to quadratures. 
In Section \ref{SecmetPP}
we shall extend the analysis to the remaining forty-four
equations with PP.

\vskip5pt We end this introduction with a comment
about the formalism used in the paper: it is elementary, and admittedly brute force (which should make the results and their proofs accessible to undergraduate students). There are other more sophisticated approaches using Cartan and tractor connections or twistor theory 
which could be 
adopted in line with \cite{BDE, East_Mat, kamran, hitchin,hitchin2}. 
\subsection*{Acknowledgements} 
We thank Phil Boalch, Robert Conte, Boris Dubrovin, Vladimir Matveev and Marta Mazzocco
for useful discussions. FC is grateful for the support of Cambridge Commonwealth, European $\&$ International Trust and CAPES Foundation Grant Proc. BEX 13656/13-9. 
MD has been partially supported by STFC consolidated grant ST/P000681/1.
\section{Proof of the main Theorem}
Our approach to proving Theorem \ref{main_theo} is based on the seminal result of Liouville
\begin{theo}[Roger Liouville 1889 \cite{Liouville1886}]
\label{theo_liouville}
A projective structure  corresponding
to the second order ODE {\em(\ref{ODE2})}
is metrisable on a neighbourhood of a point $p\in U$ iff there exist functions
$\psi_1, \psi_2, \psi_3$ defined on a neighbourhood of $p$ such that
$
\Delta\equiv\psi_1\psi_3-{\psi_2}^2\neq 0
$ 
at $p$ and the equations
\begin{subequations}
\begin{align}
\dfrac{\partial \psi_1}{\partial x}&=\frac{2}{3} A_1\psi_1-2A_0\psi_2, \label{syseq1} \\
\dfrac{\partial \psi_3}{\partial y}&=2 A_3\psi_2-\frac{2}{3}A_2\psi_3, \label{syseq2}\\
\dfrac{\partial \psi_1}{\partial y}+2\dfrac{\partial \psi_2}{\partial x}&=\frac{4}{3} A_2\psi_1-\frac{2}{3}A_1\psi_2-2A_0 \psi_3,\label{syseq3} \\
\dfrac{\partial \psi_3}{\partial x}+2\dfrac{\partial \psi_2}{\partial y}&=2 A_3\psi_1-\frac{4}{3}A_1\psi_3+\frac{2}{3}A_2 \label{syseq4}\psi_2,
\end{align}
\end{subequations}
hold on the domain of definition. The corresponding metric is then given by
\be
\label{metrisable}
g=\Delta^{-2}(\psi_1dx^2+2\psi_2dxdy+\psi_3dy^2).
\ee
\end{theo}
The system (\ref{syseq1}--\ref{syseq4})
is overdetermined, as there are more
equations than unknowns. In \cite{BDE} the integrablity
conditions were
established in terms of point invariants (\ref{ODE2}).  
The invariants obstructing metrisability vanish identically
for the projective structures arising from all six Painlev\'e equations, as these equations are non--generic
in the sense explained in \cite{BDE}: we will see that
a non--trivial solution to 
(\ref{syseq1}--\ref{syseq4}) always exists, but is 
degenerate as in general $\psi_2=\psi_3=0$. Thus the metrisability
analysis of the Painlev\'e equations needs to be carried over by  analysing the linear system 
(\ref{syseq1}--\ref{syseq4}) directly on a case by case basis.
\vskip5pt
{\bf Proof of Theorem \ref{main_theo}.}
The metrisability  
of Painlev\'e equations depend on the 
values of the parameters $(\alpha, \beta, \gamma, \delta)$. When necessary  we will indicate them in parenthesis in front of the equation label, for instance: PII$(\alpha)$, PIII$(\alpha,\beta,\gamma,\delta)$ and so on. The  Painlev\'e equations do not have a cubic term in $y'$ 
(so that $A_3=0$ in equation (\ref{ODE2}) which makes Step 1 below possible). 
A general  approach to seek solutions to the metrisability problem of this kind of projective structure is the following: 
\begin{itemize}
\item[{\bf Step 0}.] Calculate the invariants of \cite{BDE}. If they do not vanish identically, then there is no non-trivial solution to (\ref{syseq1}--\ref{syseq4}).
\item[{\bf Step 1}.] Solve equation (\ref{syseq2}) for $\psi_3$.
\item[{\bf Step 2}.] Substitute $\psi_3$ in (\ref{syseq4}) and solve it for $\psi_2$. 
\item[{\bf Step 3}.] Apply the integrability condition $\partial_x\partial_y\psi_1=\partial_y\partial_x\psi_1,\, \forall x,y,$ to the remaining equations (\ref{syseq1}) and (\ref{syseq3}).
\item[ {\bf Step 4}.] If Step 3 is successful, solve equations 
(\ref{syseq1}) and (\ref{syseq3}).
\end{itemize}
Step 0 is optional because it is equivalent to Step 3. 
After Steps 1 and 2, in general, we end up with a solution for $\psi_2$ and $\psi_3$ depending on arbitrary functions of one variable. Step 3 is then necessary to fix those functions up to constants of integration.
The above steps may be troublesome to be performed by hand, but they are easily implemented on the computer.

\vskip 5pt
We find that Painlev\'e I, II and IV are never
metrisable. On the other hand, PIII, PV and PVI are metrisable for special values of parameters, as we discuss below. The values of the parameters are found in Step 3. 
For other choices of parameters, Step 3 forces us to choose $\psi_2=\psi_3=0$ which leads to a degenerate solution. An obvious degenerate solution is the trivial one $\psi_i=0$. However, for the Painlev\'e equations, there always exist non-trivial solutions 
to the Liouville system (\ref{syseq1}--\ref{syseq4})
spanning a $1$-dimensional space, which is the maximal dimension allowed for degenerate solutions (c.f. Lemma 4.3 of \cite{BDE}). To see this, set $\psi_2=\psi_3=0$. Then (\ref{syseq1}--\ref{syseq4}) reduce to a closed overdetermined system for $\psi_1$ which has a non-vanishing solution if and only if $\partial_y A_1=2\partial_x A_2$. It is straightforward that this 
condition is fulfilled by all equations PI--PVI,
which explains 
why all invariants of \cite{BDE} vanish for Painlev\'e equations. The degenerate solutions corresponding to each Painlev\'e equation are, up to a multiplicative constant, given by
\begin{eqnarray*} 
&&
\text{PI, PII}: \psi_1=1,\quad \text{PIII}: \psi_1=\frac{y^{4/3}}{x^{2/3}}, \quad 
\text{PIV}: \psi_1= y^{2/3},\\
&&
\text{PV}: \psi_1=\frac{(1-y)^{4/3} y^{2/3}}{x^{2/3}},\quad \text{PVI}: \psi_1=(x-y)^{2/3}\left[\frac{(y-1) y}{(x-1) x}\right]^{2/3}.
\end{eqnarray*}
\begin{itemize}
\item{\bf Painlev\'e III.} Applying Steps 1 to 4 implies that a metric
exists iff 
\[
\alpha=\gamma=0\quad\mbox{or}\quad \beta=\delta=0.
\]
Both cases are essentially the same since the change of coordinates $y\mapsto {y}^{-1}$ induces PIII$(\alpha,\beta,\gamma,\delta)\to$ PIII$(-\beta,-\alpha,-\delta,-\gamma)$ and all results from one case can be recovered from the other through this map. Therefore, we only present the detailed results for $\beta=\delta=0$.
If all parameters are zero, then the projective structure is flat (which can be seen by evaluating the Liouville invariants $L_1, L_2$).
If $\beta=\delta=0$ and $(\alpha,\gamma)\not=(0,0)$
there exists  a two-dimensional family of solutions to (\ref{syseq1}--\ref{syseq4}) giving rise to the metric
\begin{eqnarray}\label{metPIII1}
g&=&
\Omega\Big(\frac{B-Axy(2\alpha+\gamma xy)}{Ax^2} dx^2+\frac{2}{xy}dx dy+
\frac{1}{y^2}dy^2
\Big),
\quad\mbox{where}\\
\Omega&=&A^{-1}(A-B+2A\alpha xy+A\gamma x^2y^2)^{-2}.\nonumber
\end{eqnarray}
The metric admits a one-parameter family of isometries 
$(x,y)\mapsto (e^\epsilon x,e^{-\epsilon} y)$. 
Setting $r=xy$ and $\theta=\ln |x|$ and rescalling the metric by $A^3$  yields
\[
g=\frac{1}{\left(-C+2\alpha r+\gamma r^2\right)^2r^2}dr^2-\frac{1}{\left(-C+2\alpha r+\gamma r^2\right)}d\theta^2,
\]
where $C=B/A-1$ is a constant. By rescaling $r$ we can set either $\alpha$ to $1$ if $\alpha\not=0$ or $\gamma$ to $\gamma/{|\gamma|}$ if $\gamma\not=0$.
\vskip5pt
If $\alpha=\beta=\gamma=\delta=0$, we have a six-dimensional family of 
solutions to (\ref{syseq1}--(\ref{syseq4}), 
all rise to  metrics of constant curvature. The projective structure is flat, and PIII$(0,0,0,0)$ can be put in the form (\ref{flat_ode}) 
with $Y=e^y$ and $X=\ln x$.
\item{\bf Painlev\'e V.} 
The projective structure is metrisable if and only if $\gamma=\delta=0$,  and is projectively 
flat if and only if $\alpha=\beta=\gamma=\delta=0$.

If $\gamma=\delta=0$ and $(\alpha,\beta)\not=(0,0)$, we have a two-dimensional family of solutions giving rise to the metric
\begin{equation}
\label{metPV1}
g=\frac{y}{A^2x^2[By+2A(\beta -\alpha y^2)]} dx^2+\frac{y}{A(y-1)^2[By+2A(\beta -\alpha y^2)]^2}dy^2,
\end{equation}
which admits $(x,y)\mapsto (e^\epsilon x, y)$ as as one-parameter family of isometries,
Defining $r=y, \theta=\ln|x|$, the metric becomes
\[
g=\frac{r}{A^3(r-1)^2[C r+2(\beta -\alpha r^2)]^2}dr^2  + 
\frac{r}{A^3[C r+2(\beta -\alpha r^2)]} d\theta^2 ,
\]
where $C=B/A$.
By redefining $C$, $B$ and $\theta$ we can set either $\beta$ to $\frac{\beta}{|\beta|}$ if $\beta\not=0$ or $\alpha$ to $\frac{\alpha}{|\alpha|}$ if $\alpha\not=0$.
\vskip5pt
If $\alpha=\beta=\gamma=\delta=0$, there exists a  six-dimensional family of solutions
to the Liouville system, each giving rise to a projectively flat metric.
Equation PV$(0,0,0,0)$ can be put in the form (\ref{flat_ode}) with $Y=\ln\left(\frac{1+\sqrt y}{\sqrt{1-y}}\right)$ and $X=\ln x$. 
\item{\bf Painlev\'e VI.}  PVI is metrisable if and only if $\alpha=\beta=\gamma=0$, $\delta=\frac{1}{2}$. In this case PVI has a solution given in terms of the elliptic integral 
\cite{Gromak1999, marta, casale07}
\begin{equation}\label{SolPVI}
\int^{y(x)}_0\frac{dw}{\sqrt{w(w-1)(w-x)}}=a\omega_1(x)+b\omega_2(x), 
\end{equation}
where the right hand side is the general solution of the Picard-Fuchs equation 
\begin{equation}\label{PFeq}
4x(x-1)\omega''(x)+4(2x-1)\omega'(x)+\omega(x)=0, 
\end{equation}
with $a$ and $b$ constants. Since the constants of integration appear linearly in (\ref{SolPVI}), the projective structure is flat\footnote{This is actually the definition of projective flatness used by Liouville \cite{Liouville1886}.}. In fact, PVI$\left(0,0,0,\frac{1}{2}\right)$ is equivalent
to (\ref{flat_ode}) in the variables $Y=\frac{1}{\omega_2(x)}\int^{y}_0\frac{dw}{\sqrt{w(w-1)(w-x)}},\; X=\frac{\omega_1(x)}{\omega_2(x)}$.
\end{itemize}
\koniec
\subsection{Coalescence}\label{SecCoalescence}
The first five Painlev\'e equations PI--PV can be derived from PVI by the process of coalescence of the parameters \cite{Ince}. In particular PIII arises from PV in the  limit $\epsilon\rightarrow 0$ 
where
\[
 x\mapsto x^2, \quad y\mapsto 1+\epsilon xy,\quad \alpha\mapsto\frac{\gamma}{8\epsilon^2}+\frac{\alpha}{4\epsilon},\quad \beta\mapsto -\frac{\gamma}{8\epsilon^2},\quad \gamma\mapsto \frac{\epsilon \beta}{4},\quad \delta\mapsto \frac{\epsilon^2\delta}{8}.
\]
We can use this process to recover the metric 
(\ref{metPIII1} ) of PIII$(\alpha,0,\gamma,0)$ from a metric of PV$(\alpha,\beta,0,0)$. To do so, it is necessary to start with (\ref{metPV1}) with the constants of integration
$$
A=\left(\frac{4 \gamma}{2 \alpha \epsilon + \gamma}\right)^{\frac{2}{3}},\quad B=\frac{(-\alpha\epsilon + \gamma)(4 \alpha \epsilon + 2 \gamma)^{\frac{1}{3}}}{
 \epsilon^2  \gamma^{\frac{1}{3}}}.
$$
Then, in the limit $\epsilon\to 0$, we find the metric (\ref{metPIII1}) with $A_{III}=1$ and $B_{III}=1-\frac{4\alpha^2}{\gamma}$, where we have attached the index $III$ to indicate that these constants $A_{III}$ and $B_{III}$ correspond to the metric of PIII$(\alpha,0,\gamma,0)$.  This is valid 
only if $\gamma\not=0$. In the case $\gamma=0$ we need $A= 4^{2/3}A_{III}$ and $B=\frac{2 \alpha A_{III} + (-A_{III} +B_{III}) \epsilon}{2^{2/3} \epsilon}$, so we still have freedom to choose two constants of integration $A_{III}$ and $B_{III}$.

\section{Reducibility and first integrals}
\label{reducibility}
The metrisable cases of PIII and PV do not define new transcendental functions, but admit a  quadrature
and are reducible to 1st order ODEs. We shall explain this in the context of Theorem \ref{main_theo} using the following  Lemma
\begin{lemma}
\label{Lemma_K}
Let 
\[
g=E(x, y)dx^2+2F(x, y)dxdy+G(x, y)dy^2
\]
be a metric on $U$ which admits a linear first integral
$K=K_1(x, y)\dot{x}+K_2(x, y)\dot{y}$.
Then
\begin{equation}\label{eq1stIntegral}
I(x, y, y')=
\frac{1}{(K_1+K_2 y')^2}\left(E+2Fy'+Gy'^2\right)
\end{equation}
is a first integral of the unparametrised geodesic equation (\ref{ODE2}).
\end{lemma}
{\bf Proof.}
Set $x^a=(x, y)$, and consider the geodesic equations for $g$ parametrised by $t$
\begin{equation}
\label{eqGeodesics}
\ddot x^a+\Gamma^a_{bc}\dot x^b \dot x^c=0.
\end{equation}
Let $t$ be a value of the affine parameter such that
$\dot{x}\neq 0$ (if no such $t$ exists then swap $x$ and $y$). Using the chain rule ${d}/{dx}=
{\dot x}^{-1}{d}/{dt}$ to eliminate $t$ between
the two equations (\ref{eqGeodesics}) yields (\ref{ODE2})
with (\ref{expressions_forA}). The geodesic Hamiltonian
$H=g_{ab}\dot{x}^a\dot{x}^b$ is a first integral
of (\ref{eqGeodesics}), but it depends on $\dot{x}$, so it does not give rise to a first integral
of (\ref{ODE2}). However dividing $H$ by the square of the linear first integral $K$ is independent on $\dot{x}^a$
and yields the first integral 
(\ref{eq1stIntegral}) for (\ref{ODE2}).
\koniec
Let us  apply this Lemma to the metrisable 
Painlev\'e cases. In case of PIII$(\alpha, 0, \gamma, 0)$ and PV$(\alpha, \beta, 0, 0)$ we shall
recover the known first integrals \cite{Gromak1999}.
\begin{itemize}
\item{\bf Painlev\'e III.} The metric (\ref{metPIII1})
admits a Killing vector
$x\p_x-y\p_y$
which gives rise to a first integral
(\ref{eq1stIntegral}) for 
PIII$(\alpha, 0, \gamma, 0)$ 
$$
I=x^2 \left(\frac{y'}{y}\right)^2+2x\frac{y'}{y}-2\alpha xy-\gamma x^2 y^2.
$$
\item{\bf Painlev\'e V.}
The metric (\ref{metPV1})
admits a Killing vector $
K=x\p_x
$
which leads to a first integral for PV$(\alpha,  \beta,  0, 0)$
\[
I=\frac{1}{y}\left(\frac{xy'}{y-1}\right)^2+\frac{2\beta}{y}-2\alpha y.
\]
\item{\bf Painlev\'e VI.}
The first integrals in this case are linear in $y'$  \cite{casale07}, and we will construct
them from the Killing vectors (rather than a quadratic integral) of the associated metric
$g=dX^2+dY^2$. The ratios of linear integrals $\dot{Y}$ and $Y\dot{X}-X\dot{Y}$ by
a linear integral $\dot{X}$ give $dY/dX$ and $Y-XdY/dX$. Evaluating these integrals
by implicitly differentiating $Y$, and using the ODE satisfied by the Wronskian
of (\ref{PFeq}) gives
\begin{equation}\label{integralPVI}
I=\frac{y'B(x)}{\sqrt{y(y-1)(y-x)}}+\int^y_0\left[A(x)+\frac{B(x)}{2(w-x)}\right]\frac{dw}{\sqrt{w(w-1)(w-x)}},
\end{equation}
where $A$ and $B$ are a solution to the Picard-Fuchs adjoint equations
\begin{align*}
\begin{cases}
A'(x)&=B(x)\frac{1}{4x(x-1)}\\
B'(x)&=-B(x)\frac{1-2x}{x(x-1)}-A(x).
\end{cases}
\end{align*}
\end{itemize}

 A prolongation of the metrisability equations 
(\ref{syseq1}-\ref{syseq4}) leads to a closed system of six linear PDEs for six unknowns \cite{BDE}. The dimension $m([\nabla])$ of the vector space of solutions 
to this system is called the degree of mobility of the projective structure.
In the generic, non--metrisable case $m([\nabla])=0$, and in the projectively flat case
$m([\nabla])=6$. The Koenigs theorem \cite{Koenigs} states that $m([\nabla])\neq 5$.
The construction below applies to projective structures
where  $m([\nabla])>1$.
\begin{prop}\label{thm1stintegral}
If a projective structure $[\nabla]$ in two  
dimensions admits  two linearly independent 
solutions $\psi^{(1)}$ and $\psi^{(2)}$ to the metrisability equations 
(\ref{syseq1}-\ref{syseq4}), then
\be
\label{quadraticI}
I(x,y,y'):=\frac{\psi^{(1)}_{1}+2\psi^{(1)}_{2}y'+\psi^{(1)}_{3}y'^2}{\psi^{(2)}_{1}+2\psi^{(2)}_{2}y'+\psi^{(2)}_{3}y'^2}
\ee
is a first integral of the unparametrised geodesic equation (\ref{ODE2}).

If there  exists a linear combination  of $\psi^{(1)}$ and $\psi^{(2)}$ which is 
degenerate, 
then any metric $g$ compatible with $[\nabla]$ admits a Killing vector.
\end{prop}
{\bf Proof.} The constancy of (\ref{quadraticI}) could be established by explicitly evaluating
$dI$ on solutions to (\ref{syseq1}-\ref{syseq4}), which gives $0$. Below, we shall use a less direct method
which will allow us to prove both parts of the Proposition.
Two connections
$\nabla$ and $\hat{\nabla}$ belong to the same projective equivalence class 
$[\nabla]$ if there exists a 
one-form $\Upsilon$ on $U$ such that (\ref{projequiv}) holds.
Consider a connection\footnote{ 
A more invariant way to define this connection is as follows. 
Pick a connection $\nabla\in[\nabla]$ and set 
\[
\Pi^c_{ab}=\Gamma^c_{ab}-\frac{1}{3}\Gamma^d_{da}\delta^c_b-\frac{1}{3}\Gamma^d_{db}\delta^c_a.
\]
The object $\Gamma^d_{da}$ is not a 1-form, and thus $\Pi^c_{ab}$ does not transform as an affine connection in general, but only under coordinate transformations of constant Jacobian. So once we choose this representative we can only apply this kind of coordinate 
transformations in (\ref{eqmetsym}). Thomas \cite{Thomas1925} introduced the terminology {\em equi-transformation} 
for coordinate 
changes preserving the volume (of Jacobian identically $1$), 
{\em projective connection} for $\Pi^c_{ab}$ and 
{\em equi-tensor} for entities such as $\Gamma^d_{ad}$ transforming 
like tensors under equi-transformations.}
$D\in[\nabla]$ with Christoffel symbols given by
\[
\Pi_{11}^1=\frac{1}{3}A_1, \quad \Pi_{12}^1=\frac{1}{3}A_2, \quad \Pi_{22}^1=
A_3, \quad \Pi_{11}^2=-A_0, \quad \Pi_{21}^2=-\frac{1}{3}A_1, \quad
\Pi_{22}^2=-\frac{1}{3}A_2.
\]
Set
$\psi_1=\sigma_{11}$, $\psi_2=\sigma_{12}$ 
and $\psi_3=\sigma_{22}$. 
Then the metrisability equations (\ref{syseq1}-\ref{syseq4})
are equivalent to the Killing tensor equation \cite{BDE}
\begin{equation}\label{eqmetsym}
D_{(a}\sigma_{bc)}=0.
\end{equation}
Therefore $\sigma^{(i)}$ $(i=1,2)$ are Killing tensors, 
and $I^{(i)}(t):=\sigma^{(i)}_{11}\dot x^2+2\sigma^{(i)}_{12}\dot x\dot y+\sigma^{(i)}_{22}\dot y^2$ are conserved along geodesics of $D$
$$
\dfrac{d}{dx}I(x,y(x),y'(x))=\frac{1}{\dot x}\dfrac{d}{dt}\frac{I^{(1)}(t)}{I^{(2)}(t)}=0,
$$
where we have used $\dot y/\dot x=y'$ to write $I=I^{(1)}/I^{(2)}$.

For the second part, the projective structure is metrisable 
(this is true even if both $\psi^{(i)}$ are degenerate, as there always exists a non--degenerate
linear combination, i. e. two degenerate solutions can only differ by a constant multiple. See Lemma 4.3 in \cite{BDE}). Without loss of generality say
that $\psi^{(2)}$ is degenerate. Then there exists a non-vanishing one-form 
$\omega$ such that $\sigma^{(2)}_{ab}=\omega_a\omega_b$. Then the metrisability equations 
(\ref{eqmetsym}) yield
$$
D_{(a}\omega_{b)}=0.
$$
The Levi-Civita connection $\nabla$ of the metric 
$g=\sigma/\Delta^2$
is obtained from $D$ by applying a transformation (\ref{projequiv}) with an equi-tensor 
$\Upsilon_a=\nabla_a \left(-\frac{1}{2}\ln\left|\Delta\right|\right)$, 
and we verify that
\[
\nabla_{(a}K_{b)}=0,
\]
where $K={\Delta}^{-1}\omega$. Thus $K$ is a linear first integral of the geodesic flow of $g$.
\koniec 
{\bf Remarks.}
\begin{itemize}
\item
Not all projective structures with $m([\nabla])>1$ admit a linear first integral. The metrics 
\cite{dini}
$$
g_1=\left(X(x)-Y(y)\right)\left(dx^2+dy^2\right) \quad \text{ and } \quad g_2=\left(\frac{1}{Y(y)}-\frac{1}{X(x)}\right)\left(\frac{dx^2}{X(x)}+\frac{dy^2}{Y(y)}\right),
$$
are projectively equivalent with an unparametrised geodesic equation
\begin{equation}\label{eqnoKilling}
y''+\frac{1}{2\left(X(x)-Y(y)\right)}\left(Y'+X' y'+Y' y'^2+X' y'^3\right)=0.
\end{equation}
These metrics in general do not admit a Killing vector, but clearly $m([\nabla])>1$. The 
first
 integral (\ref{quadraticI}) is
$$
I=\frac{Y(y)+X(x)y'^2}{1+y'^2}.
$$
\item
Each Painlev\'e equation admits a degenerate solution
to the metrisability equations. This implies
that the corresponding projective class $[\nabla]$ contains
a representative $\nabla$ which has  symmetric Ricci tensor, and
admits a linear first integral. 
In \cite{FCMD2015} it was shown that
for such affine connections $\nu_5=0$,  where $\nu_5$ is
a point invariant for (\ref{ODE2}) defined by Liouville
\cite{Liouville1886}. This is in agreement with 
\cite{HietDry2002}, where it was stated that $\nu_5$ vanishes for all Painlev\'e equations.
\item In \cite{mettler} it was shown that all two--dimensional projective structures
are locally Weyl--metrisable. For a given ODE (\ref{ODE2}) finding an explicit expression for the Weyl connection reduces to constructing a point transformation such that
$A_0=A_2$ and $A_1=A_3$. This should in principle be possible
of all six Painleve equations, but the resulting ODEs
may not have Painleve property if the point transformation in not
of the form (\ref{homot}).
\item  In the recent work \cite{Levi} some connections between
the Painlev\'e property and Lie point symmetries have been uncovered.
While the problems studied in \cite{Levi} are different than those addressed in our work,
some of the results appear to be related. In particular among the six Painlev\'e transcendents only  PIII and PV have nontrivial symmetry algebras and that only for special values of the parameters.
\end{itemize}

\section{Metrisability of equations with the 
Painlev\'e property}\label{SecmetPP}

All fifty equivalence classes of 2nd order ODEs with 
Painlev\'e property are
of the form (\ref{ODE2}) and so they define projective 
structures. Six of them are the Painlev\'e equations
and their metrisability is determined by Theorem
\ref{main_theo}.
In this Section we summarise the results of the analysis of the remaining forty-four cases listed in 
\cite{Ince} in their most general form. We use the same numbering as this reference. We can divide these equations in five sets, according to their metrisability properties:
\begin{itemize}
\item[1.] Metrisable with one degenerate solution and $4>m([\nabla])>1$: II, III, VII, VIII, XII, XVIII, XIX, XXI, XXIII, XXIX, XXX, XXXIII, XXXVIII, XLIV, XLIX.
\item[2.] Metrisable with degenerate solution and  $m([\nabla])=4$ : XXII, 
XXXII.
\item[3.] Non-metrisable, but  admitting a degenerate solution: XIV, XX, XXXIV.
\item[4.] Not metrisable, and no non-trivial 
solutions to the metrisability equations: V, X, XV, XVI, XXIV, XXV, XXVI, XXVII, XXVIII, XXXV, XXXVI, XL, XLII, XLVII, XLV, XLVI, XLVIII.
\item[5.] Metrisable and projectively flat: I, VI, XI, XVII, XXXVII, XLI, XLIII.
\end{itemize}
The Painlev\'e equations are IV, IX, XIII, XXXI, XXXIX, L, which we did not include 
in the list but would fit in group 3 in general. 
The metrisable cases all admit a degenerate solution (thus their metrics admit a 
Killing vector, from Proposition  \ref{thm1stintegral}), and their ODEs admit a quadratic first integral. 

The submaximal (i.e. degree of mobility  $m([\nabla])=4$)
equations XXII and XXXII are related by a point transformation which is however not rational.  
The ODE XXXII 
\[
y''=\frac{1}{2y}(1+(y')^2)
\]
is metrisable by
$g=y(dx^2+dy^2)$. The four quadratic first integrals for the parametrised geodesic motion
give rise to three functionally dependent integrals quadratic in $y'$. Two independent integrals are
\[
I_1=\frac{1}{y}(1+(y')^2), \quad I_2=2y'-\frac{x}{y}(1+(y')^2).
\]
\section{Summary}
We have established which 2nd order ODEs with Painlev\'e
property are metrisable, i. e. all their integral curves
are geodesics of some (pseudo) Riemannian metric.
Out of the six Painlev\'e equations
only 
PIII$(\alpha, 0, \gamma, 0)$, 
PIII$(0, \beta, 0, \delta)$, 
PV$(\alpha, \beta, 0, 0)$ and
PVI$(0, 0, 0, 1/2)$ are metrisable, the last case being
projectively flat. In all cases the metrisable equations
with PP admit a first integral, and
the degree of mobility is at least two. Thus metrisability picks out non--transcendental cases
in the Painlev\'e analysis.

It would be interesting to extend Theorem \ref{main_theo} to systems of two second order ODEs
\be
\label{systems}
y''=F(x, y, z,  y', z'), \quad z''=G(x, y, z, y', z').
\ee
It is known how to characterise the systems resulting from a three--dimensional projective 
structure \cite{Fels, CDT13, Doubrov}, and some necessary and sufficient conditions for metrisability have recently 
been constructed \cite{DE} and \cite{E}. The classification of systems (\ref{systems})
which admit Painlev\'e property is however missing.


\begin{thebibliography}{99}

\bibitem{BDE} Bryant, R. L., Dunajski, M., and Eastwood, M. G. (2009) 
Metrisability of two-dimensional projective structures,
J. Differential Geometry {\bf 83}, 465--499.







\bibitem{C22} Cartan, E. (1924) 
Sur les vari\'et\'es \`a connexion projective,
Bull. Soc. Math. France {\bf 52}, 205--241.

\bibitem{casale07} Casale, G. (2007)
The Galois groupoid of Picard-Painlev{\'e} VI equation,
Proceedings of the French-Japanese Conference Algebraic, Analytic and Geometric Aspect of Complex Differential Equations and their Deformations. 
Painlev{\'e} Hierarchies.
Research Institute for Mathematical Sciences.

\bibitem{CDT13} Casey, S., Dunajski, M. and Tod, K. P. (2013)
Twistor geometry of a pair of second order ODEs. Comm. Math. Phys. 321, 681–701
  
\bibitem{FCMD2015} Contatto, F. and Dunajski, M. (2016).
First integrals of affine connections and Hamiltonian systems of 
hydrodynamic type.  Journal of Integrable Systems {\bf 1.}
{\tt arXiv:1510.01906}.
	



\bibitem{dini} Dini, U. (1869) 
Sopra un problema che si presenta nella teoria generale delle rappresentazioni geografiche di una superficie su un altra.
Ann. Math. Ser 2. {\bf 3} 269-293.

\bibitem{Doubrov}
Doubrov, B., \& Medvedev, A.  (2014)
Fundamental invariants of systems of ODEs of higher order
Diff. Geom. App. {\bf 35}, 291–313.
\bibitem{DE} Dunajski, M. \& Eastwood, M. G. (2016)
Metrisability of three-dimensional path geometries.
European Journal of Mathematics. {\tt arXiv:1408.2170}

\bibitem{East_Mat}  Eastwood, M. G. \&  Matveev, V (2007)
Metric connections in projective differential geometry, in
\emph{Symmetries and Overdetermined Systems of Partial Differential Equations},
IMA Volumes in Mathematics and its Applications 144, 
Springer Verlag 2007, pp.~339--350.

\bibitem{E} Eastwood, M. (2017) Metrisability of three-dimensional projective structures. {\tt  arXiv:1712.06191}.

\bibitem{Fels} Fels, M. (1995)
The equivalence problem for systems of second-order ordinary differential equations, Proc. London Math. Soc.  {\bf 71} , 221-240.


\bibitem{Gromak1999} Gromak, V. I. (1999)
B\"acklund Transformations of Painlev\'e Equations and Their Applications.
The Painlev\'e Property, One Century Later, Conte, R.
CRM Series in Mathematical Physics

\bibitem{HietDry2002} Hietarinta, J. \& Dryuma, V. (2002) 
Is my ODE a Painlev\'e equation in disguise? 
J. Nonlinear Math. Phys. {\bf 9},  suppl. 1, 67--74.


\bibitem{hitchin} Hitchin, N. J. (1997) 
Geometrical aspects of Schlesinger's equation. J. Geom. Phys. {\bf 23}

\bibitem{hitchin2} Hitchin, N.~J. (1982)
Complex manifolds and {E}instein's equations.
In {\em Twistor geometry and Non-Linear Systems\/}, Lecture
Notes in Math., 970, Springer, pp.~73--99.


\bibitem{Ince} {Ince, E. L.} (1965)
{\it Ordinary Differential Equations}, Dover.
  
\bibitem{kamran}
Kamran, N. Lamb, K. G. Shadwick, W. F. (1985)
 The local equivalence problem for 
$d^2y/dx^2=F(x,y,dy/dx)$ and the Painlev\'e transcendents. 
J. Differential Geom. {\bf 22}, 139–150. 

\bibitem{kartak} Kartak, V. V. (2013)
 Point classification of second order ODEs and its application to Painlevé equations. J. Nonlinear Math. Phys. {\bf 20}.


\bibitem{Koenigs} Koenigs, M. G.
{Sur les g\'eodesiques a int\'egrales quadratiques.}
Note II from Darboux' 
{\em Le\c{c}ons sur la th\'eorie g\'en\'erale des surfaces}, 
 Vol. IV, Chelsea Publishing, 1896.

\bibitem{Levi}
Levi, D. Sekera, D., and Winternitz, P. (2017)
Lie point symmetries and ODEs passing the Painlev\'e test.
{\tt arXiv:1712.09811}

\bibitem{Liouville1886} Liouville, R. (1889)
Sur les invariants de certaines \'equations diff\'erentielles et
sur leurs
applications,
Jour. de l'Ecole Politechnique, {\bf Cah.59}, 7--76.

\bibitem{marta} Mazzocco, M. (2001)
Picard and Chazy solutions to the Painlev\'e VI equation.
Math. Ann.  {\bf 321}, 157-195.


\bibitem{mettler} Mettler, T. (2014)
Weyl metrisability of two-dimensional projective structures.
Math. Proc. Cambridge Philos. Soc. {\bf 156}, 99-113

\bibitem{painleve} Painlev\'e, P. (1902)
Sur les équations différentielles du second ordre et d'ordre superieur dont 
l'integrale generale est uniformee,  Acta Math. {\bf 25}: 1–85.


\bibitem{Thomas1925}
Thomas, T. Y. (1925) 
On the Projective and Equi-Projective Geometries of Paths,
Proceedings of the National Academy of Sciences
{\bf 11}.

\end{thebibliography}
\end{document}